\documentclass[10pt]{tran-l}
\usepackage{graphicx,amscd,amssymb,verbatim}
\usepackage[dvips]{hyperref}
\usepackage[leqno]{amsmath}
\usepackage{subfigure}

 \newtheorem{thm}{Theorem}
 \newtheorem{Corollary}{Corollary}
 \newtheorem{Lemma}{Lemma}

 \newcommand{\Real}{\mathbb{R}}

\begin{document}

\title[One group of inequalities]{One group of inequalities with altitudes and medians in triangle}
\author{Zhivko Zhelev}

\address{Zhivko Zhelev\newline
\hspace*{5mm}Department of Mathematics and Informatics,\newline
 \hspace*{5mm}Trakia University, Rektorat, AF, 367A,\newline
  \hspace*{5mm}6003 St. Zagora,\newline
  \hspace*{5mm}Bulgaria.\newline
   \hspace*{5mm}Email: \tt{zhelev@uni-sz.bg}}

\date{}

\begin{abstract}
In the article we prove some inequalities that contain relations
between altitudes and medians in triangle. At least one of these
inequalities has not been considered in the literature before and
the main theorem has also not been proved elsewhere in that
form.\\
\indent Some immediate corollaries have been presented as well.
\end{abstract}

\keywords{medians, altitudes, inequality.
\newline \hspace*{5mm}2000 {\it Mathematics Subject Classification.} 51M04, 51M05, 51M16.}

 \maketitle
\indent Geometry of the triangle is a realm of elementary geometry
where interested new results pop up all the time. There are plenty
of theorems concerning geometry of the triangle, including
hundreds of geometric inequalities (see for example
\cite{Bottema}). The following result consists of two inequalities
which look quite pleasant but the second one turned out to be
extremely difficult to tackle. In fact we prove the following
\begin{thm}
Let $\triangle ABC$ be an arbitrary triangle with sides $a$, $b$
and $c$. Let also $h_a, h_b, h_c$ and $m_a, m_b, m_c$ be the
altitudes and medians to these sides respectively. Then the
following two inequalities hold:

\begin{eqnarray}
 ah_a+bh_b+ch_c\le\sqrt{bc}h_a+\sqrt{ac}h_b+\sqrt{ab}h_c\\
  am_a+bm_b+cm_c\le\sqrt{bc}m_a+\sqrt{ac}m_b+\sqrt{ab}m_c .
  \end{eqnarray}
\end{thm}

\indent As far as the author knows, the second inequality has not
been proved yet. But different random checks of values of the
triangle sides have not brought any counterexamples.\\
\indent Our proof is based on  one and two variable functions
theory which we consider an impediment since the problem
formulated above is in the area of the elementary mathematics.\\
\indent First we formulate some lemmas which are more or less
obvious.

\begin{Lemma}[{\bf{Arithmetic-Geometric-Mean Inequality}}] Let $a_1,
a_2,\ldots a_n\in\Real$ and $a_i\ge 0$, $i=1,\ldots ,n$. Then the
following inequality is true:
\begin{equation}
\dfrac{a_1+a_2+\cdots + a_n}{n}\ge \root n\of{a_1a_1\cdots a_n} .
\end{equation}
\end{Lemma}

\medskip

\underline{\bf{Proof.}} See [4, p. 18-19] and \cite{Nelson} for a
detailed proof.$\hfill\square$

\begin{Lemma}
If $a+b+c>0$, then $a^3+b^3+c^3\ge 3abc$.
\end{Lemma}

\underline{\bf{Proof.}} From Lemma 1 ($n=3$), we get that
$a_1+a_2+a_3\ge 3\root 3\of{a_1a_2a_3}$. Replacing $a_1$ with
$a^3$, $a_2$ with $b^3$ and $a_3$ with $c^3$, we get that
$a^3+b^3+c^3\ge 3\root 3\of{a^3b^3c^3}=3abc$., q. e.
d.$\hfill\square$\\

\bigskip
\underline{\bf{Proof of the main theorem.}} \\

In order to prove (1) we use that $ah_a=bh_b=ch_c=2S_{\triangle
ABC}$. Then
$$
\begin{array}{c}
ah_a+bh_b+ch_c-\sqrt{bc}h_a-\sqrt{ac}h_b-\sqrt{ab}h_c=
6S-\dfrac{2S\sqrt{bc}}{a}-\dfrac{2S\sqrt{ac}}{b}-\dfrac{2S\sqrt{ab}}{c}=\\
\\
2S\left
(3-\dfrac{\sqrt{bc}}{a}-\dfrac{\sqrt{ac}}{b}-\dfrac{\sqrt{ab}}{c}\right
).
\end{array}
$$
Therefore,
$$
\begin{array}{c}
ah_a+bh_b+ch_c\le\sqrt{bc}h_a+\sqrt{ac}h_b+\sqrt{ab}h_c\iff\left
(3-\dfrac{\sqrt{bc}}{a}-\dfrac{\sqrt{ac}}{b}-\dfrac{\sqrt{ab}}{c}\right
)\le 0\\
\\
\iff 3abc\le
(bc)^{\frac{3}{2}}+(ac)^{\frac{3}{2}}+(ab)^{\frac{3}{2}},\,\,a>0,\,\,
b>0,\,\, c>0.
\end{array}
$$

\indent By the substitution $x=\sqrt{bc}$, $y=\sqrt{ac}$,
$z=\sqrt{ab}$, we get that
$$
3abc\le
(bc)^{\frac{3}{2}}+(ac)^{\frac{3}{2}}+(ab)^{\frac{3}{2}}\iff
3xyz\le x^3+y^3+z^3,\,\,x>0,\,\, y>0,\,\, z>0.
$$
\medskip

\indent Last inequality is exactly Lemma 2 and the proof of (1) is
completed.\\
\indent Now, without loss of generality, we can assume for
$\triangle ABC$, that $a\ge b\ge c$. Three main cases are
possible:
\begin{itemize}
\item[$\dagger )$] $\triangle ABC$ is equilateral, i. e. $a=b=c$;
\item[$\dagger\dagger )$] $\triangle ABC$ is isosceles, i. e.
$a=b>c$; \item[$\dagger\dagger\dagger )$] $\triangle ABC$ is an
arbitrary triangle, i. e. $a>b>c$.
\end{itemize}
\medskip

\indent First, let's rewrite (2) in the form
\begin{equation}
(a-\sqrt{bc})m_a+(b-\sqrt{ac})m_b+(c-\sqrt{ab})m_c\le 0.
\end{equation}

\medskip

\underline{\bf{First case.}} If $\triangle ABC$ is equilateral,
then $a-\sqrt{bc}=b-\sqrt{ac}=c-\sqrt{ab}=0$ and (4) is fulfilled
as an equality. We will see below that this case is in fact the
extremal case for our problem.\\

\bigskip

\underline{\bf{Second case.}}Now let $\triangle ABC$ be isosceles
and we may assume $a=b>c>0$. On the other hand we have that

\begin{equation}
\begin{array}{c}
m_a=\dfrac{1}{2}\sqrt{2b^2+2c^2-a^2}=\dfrac{1}{2}\sqrt{a^2+2c^2}=m_b,\\
m_c=\dfrac{1}{2}\sqrt{2a^2+2b^2-c^2}=\dfrac{1}{2}\sqrt{4a^2-c^2}.
\end{array}
\end{equation}

\indent Then using (5), (4) is transformed into

\begin{equation}
(a-\sqrt{ac})\sqrt{a^2+2c^2}\le\dfrac{a-c}{2}\sqrt{4a^2-c^2},\quad
a>c.
\end{equation}

We will prove (6). After some tedious computations, we get
consequently:

$$
\begin{array}{c}
(a-\sqrt{ac})\sqrt{a^2+2c^2}\le\dfrac{a-c}{2}\sqrt{4a^2-c^2}\iff\\
\\
\sqrt{a}(\sqrt{a}-\sqrt{c})\sqrt{a^2+2c^2}\le\dfrac{(\sqrt{a}-\sqrt{c})(\sqrt{a}+\sqrt{c})}{2}
\sqrt{4a^2-c^2}\iff\\
\\
2\sqrt{a^3+2ac^2}\le (\sqrt{a}+\sqrt{c})\sqrt{4a^2-c^2}\iff\\
\\
\dfrac{4(a^3+2ac^2)}{4a^2-c^2}\le
a+c+2\sqrt{ac}\iff\dfrac{4a^3+8ac^2-4a^3+ac^2-4a^2c+c^3}{4a^2-c^2}\le
2\sqrt{ac}\iff\\
\\
\left (\dfrac{c^3+9ac^2-4a^2c}{4a^2-c^2}\right )^2\le 4ac\iff\\
\\
c^6+81a^2c^4+16a^4c^2+18ac^5-8a^2c^4-72a^3c^3\le
64a^5c+4ac^5-32a^3c^3\iff\\
c^6+14ac^5+73a^2c^4-40a^3c^3+16a^4c^2-64a^5c\le 0\iff\\
c(c^5+14ac^4+73a^2c^3-40a^3c^2+16a^4c-64a^5)\le 0\iff\\
c^5+14ac^4+73a^2c^3-40a^3c^2+16a^4c-64a^5\le 0\iff\\
\\
\left (\dfrac{c}{a}\right )^5+14\left (\dfrac{c}{a}\right )^4+
73\left (\dfrac{c}{a}\right )^3-40\left (\dfrac{c}{a}\right
)^2+16\left (\dfrac{c}{a}\right )-64\le 0.
\end{array}
$$

\indent In the last expression, let $\dfrac{c}{a}=t\in (0,1)$, and
it follows that

$$
\begin{array}{c}
t^5+14t^4+73t^3-40t^2+16t-64\le 0\iff\\
\\
\underbrace{(t-1)}_{<0}\underbrace{(t^4+15t^3+88t^2+48t+64)}_{>0\,\,\,\mbox{if}\,\,\,
t>0}\le 0.
\end{array}
$$

But the last one is obviously true and that proves (6) in this
case.$\hfill\square$\\

\bigskip
\underline{\bf{Third case.}} Let now $\triangle ABC$ be an
arbitrary triangle and let $a>b>c>0$. We rewrite (4) in the form

\begin{equation}
\dfrac{1}{2}(a-\sqrt{bc})\sqrt{2b^2+2c^2-a^2}+\dfrac{1}{2}(b-\sqrt{ac})
\sqrt{2a^2+2c^2-b^2}+\dfrac{1}{2}(c-\sqrt{ab})\sqrt{2a^2+2b^2-c^2}\le
0
\end{equation}

or

\begin{equation}
\begin{array}{c}
\dfrac{1}{2}a^2\underbrace{\left [\left
(1-\sqrt{\dfrac{b}{a}}\sqrt{\dfrac{c}{a}}\right )\sqrt{2\left
(\dfrac{b}{a}\right )^2+2\left (\dfrac{c}{a}\right )^2-1}+ \left
(\dfrac{b}{a}-\sqrt{\dfrac{c}{a}}\right )\sqrt{2+2\left
(\dfrac{c}{a}\right )^2-\left (\dfrac{b}{a}\right )^2}\right ]}_
{F\left (\dfrac{b}{a},\dfrac{c}{a}\right ):=F(x,y)}+\\
 \underbrace{\dfrac{1}{2}a^2\left [\left
(\dfrac{c}{a}-\sqrt{\dfrac{b}{a}}\right )\sqrt{2+2\left
(\dfrac{b}{a}\right )^2-\left (\dfrac{c}{a}\right )^2}\right
]}_{F\left (\dfrac{b}{a},\dfrac{c}{a}\right ):=F(x,y)}\le 0,
\end{array}
\end{equation}

\noindent and therefore $\dfrac{1}{2}a^2F(x,y)\le 0\iff F(x,y)\le 0$.\\
\indent Two variable function $F(x,y)$ defined above, we name
\textit{devil-fish function} and the surface this function plots
in $\Real^3$ -- \textit{devil-fish surface}.\\
\indent In our new notation inequality (4) transforms into the
extremal problem
\begin{equation}
\begin{array}{c}
\displaystyle{\max_{(x,y)\in M\subset\mathbb{R}^2}F(x,y)},\mbox{where}\\
\\
M=\{(x,y)\in\mathbb{R}^2|\,\,0\le y\le x\le 1,\,\,x+y\ge 1\},
\end{array}
\end{equation}

\noindent where $M$ is geometrically a right-angle triangle and
the last inequality follows from the fact that $a+b>c$ in an
arbitrary
triangle.\\
\indent It is a straightforward check that the devil-fish function
is well-defined on $M$. Therefore, in order to prove (4), we have
to prove that $\displaystyle{\max_{(x,y)\in M}F(x,y)}=0$.\\

\medskip
Note that $M$ is a compact set and since $M$ is a continuous
function on $M$, it follows from some classical results in the
real mathematical analysis, that $F(x,y)$ reaches its maximum
there. In order to find that maximal point we consider the
interior and the boundary of $M$ separately, $M=\mbox{int}\,
M\cup\partial
M$.\\

\indent Using some calculating programs such as \textit{Maple}
(computations can be done manually but that can turn into an
extremely tedious hardwork), we find that

$$
\begin{array}{c}
\dfrac{\partial F(x,y)}{\partial
x}=-\dfrac{y\sqrt{2x^2+2y^2-1}}{2\sqrt{xy}}+\dfrac{2(1-\sqrt{xy})x}{\sqrt{2x^2+2y^2-1}}+\sqrt{2+2y^2-x^2}-\\
\\
\dfrac{(x-\sqrt{y})x}{2+2y^2-x^2}-\dfrac{2+2x^2-y^2}{2\sqrt{x}}+\dfrac{2(y-\sqrt{x})x}{\sqrt{2+2x^2-y^2}},
\end{array}
$$

\noindent and since the devil-fish function is a symmetric
function, i. e. $F(x,y)=F(y,x)$, it follows immediately that

$$
\begin{array}{c}
\dfrac{\partial F(x,y)}{\partial
y}=-\dfrac{x\sqrt{2y^2+2x^2-1}}{2\sqrt{xy}}+\dfrac{2(1-\sqrt{xy})y}{\sqrt{2y^2+2x^2-1}}+\sqrt{2+2x^2-y^2}-\\
\\
\dfrac{(y-\sqrt{x})y}{2+2x^2-y^2}-\dfrac{2+2y^2-x^2}{2\sqrt{y}}+\dfrac{2(x-\sqrt{y})y}{\sqrt{2+2y^2-x^2}}.
\end{array}
$$

\indent Then we solve the system

$$\left|
\begin{array}{c}\dfrac{\partial F}{\partial x}=0\\
\\
 \dfrac{\partial F}{\partial y}=0\end{array}, \right.
 $$
\medskip

 \noindent which gives us two points: $M_1(0,9238127491\ldots ,0,1660179102\ldots )$
 and
 $M_2(1,1)$. Note that $M_1, M_2\in M$ and $M_1\in\mbox{int}\,M$, à $M_2\in\partial M$.

\indent Again using some software, we get for the devil-fish
function's hessian:

$$
 H(x,y):=\dfrac{\partial^2 F}{\partial x^2}\cdot\dfrac{\partial^2
 F}{\partial y^2}-\left [\dfrac{\partial^2F}{\partial x\partial y}\right
 ]^2_{|M_2(1,1)}=\left (-\dfrac{3\sqrt{3}}{2}\right )\left
 (-\dfrac{3\sqrt{3}}{2}\right )-\left (\dfrac{3\sqrt{3}}{4}\right
 )^2=\left (\dfrac{9}{4}\right )^2>0
 $$

 \noindent and since $\dfrac{\partial^2F}{\partial
  x^2}_{|M_2(1,1)}=-\dfrac{3\sqrt{3}}{2}<0$, it follows that point
  $M_2(1,1)$ is a possible point of maximum over $M$. For
  completeness we need to check what is going on the boundary
  because $M_2\in\partial M$. But that is straightforward and this
  checking consists of the cases $\{x=0\}$, $\{x=y\}$, $\{x=1\}$,
  and $\{x+y=1\}$. This leads to the fact that the point in
  question is in fact an absolute maximum in $M$. This yields

  $$
  \sup_{(x,y)\in M} F(x,y)=\max_{(x,y)\in M} F(x,y)=F(1,1)=0.
  $$
\medskip

Additionally, $M_1$ is a point of a local minimum and moreover:

$$
  F(x,y)_{|M_1}=F(0,9238127491\ldots ,0,1660179102\ldots
  )=-0,4280657968\ldots
  $$
  \newpage
  The devil-fish surface can be seen below:

  \begin{figure}[ht]
     \centering
       \includegraphics [keepaspectratio, scale=0.5]{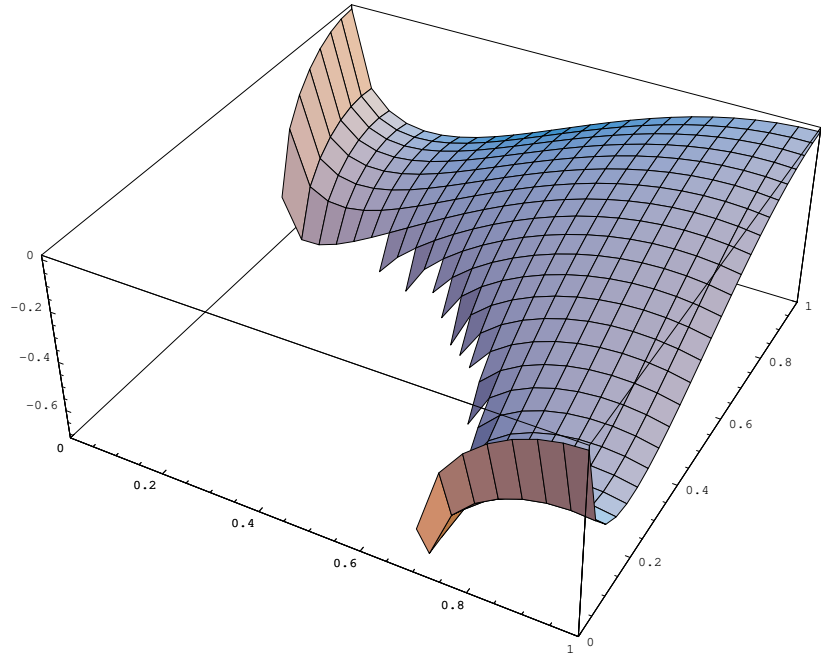}
  \caption{Devil-fish surface}
  \end{figure}

  Using software such as \textit{Maple} and \textit{Mathematica},
  one can get different views to that surface:
\newpage
  \begin{figure}
   \centering
     \subfigure[View from above]{
     \includegraphics[keepaspectratio,scale=0.35]{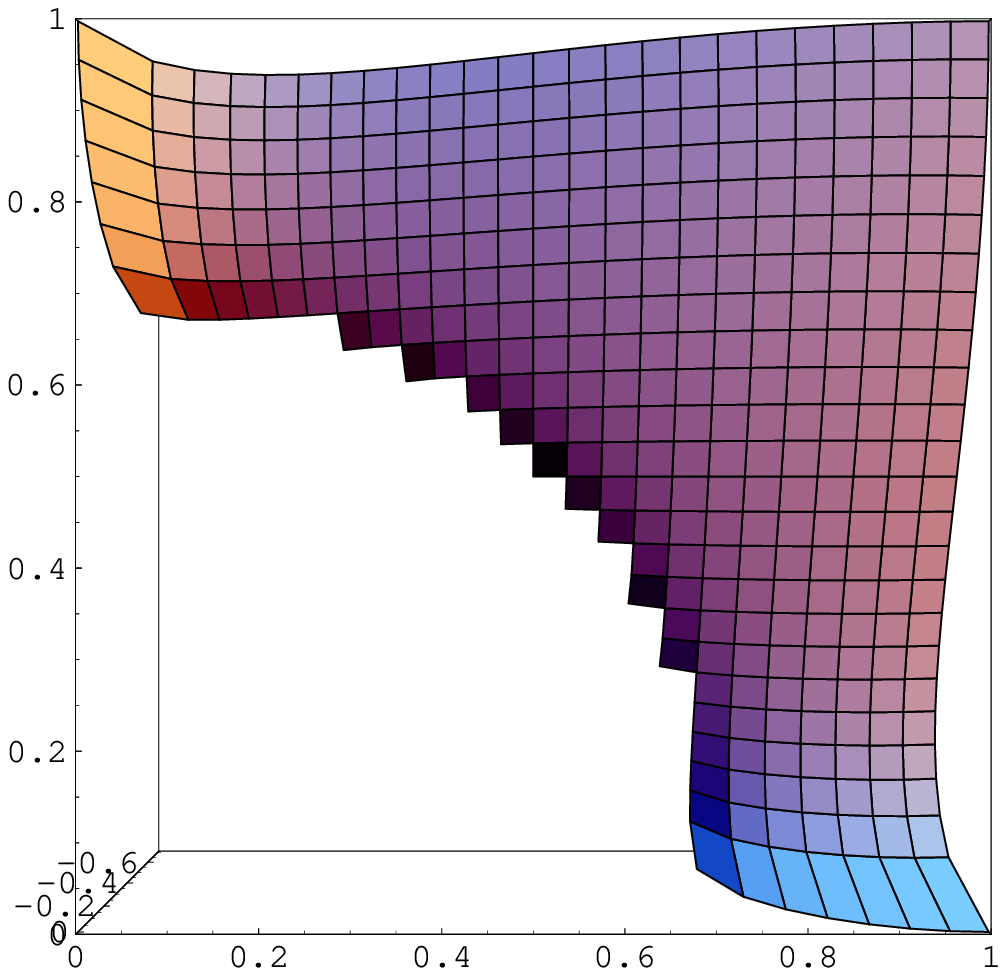}}
\subfigure[Left corner view]{
   \includegraphics[keepaspectratio,scale=0.35]{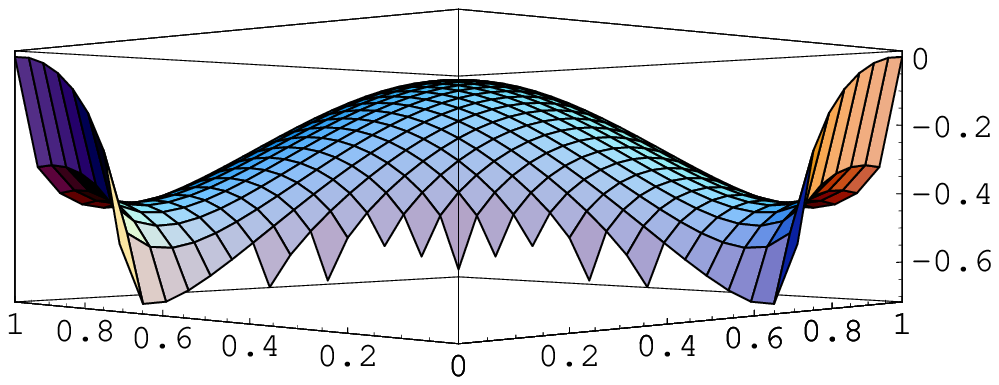}}
\subfigure[Right corner view]{
   \includegraphics[keepaspectratio,scale=0.35]{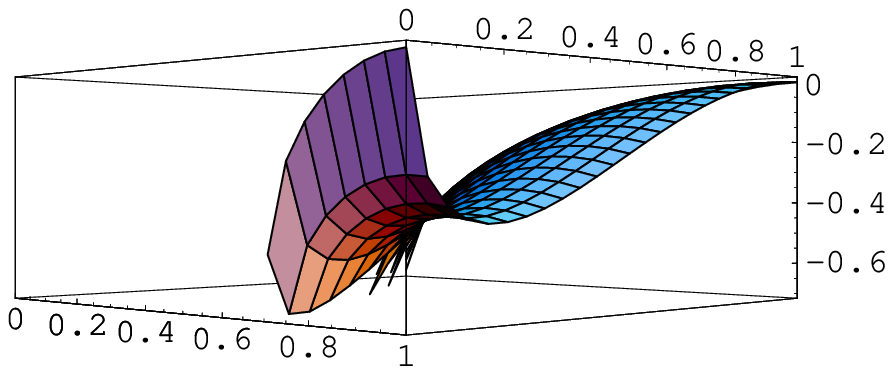}}
   \end{figure}

\begin{figure}
  \centering
   \subfigure[Front view]{
   \includegraphics[keepaspectratio,scale=0.35]{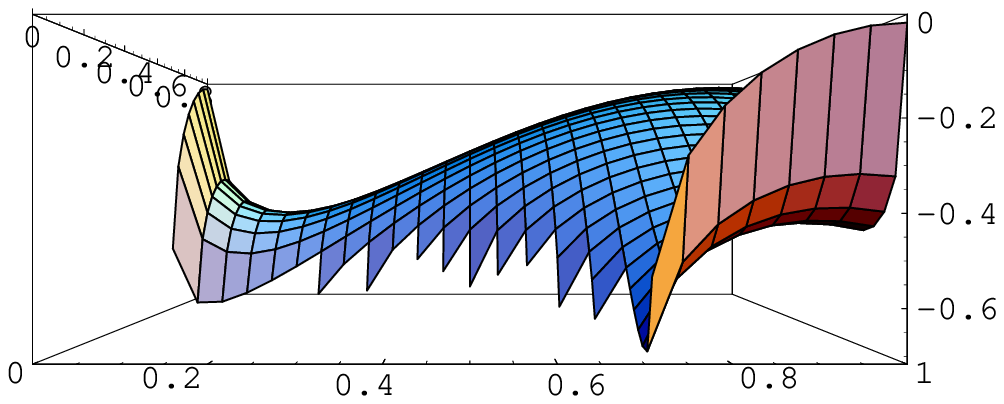}}
\subfigure[Front-down view]{
   \includegraphics[keepaspectratio,scale=0.35]{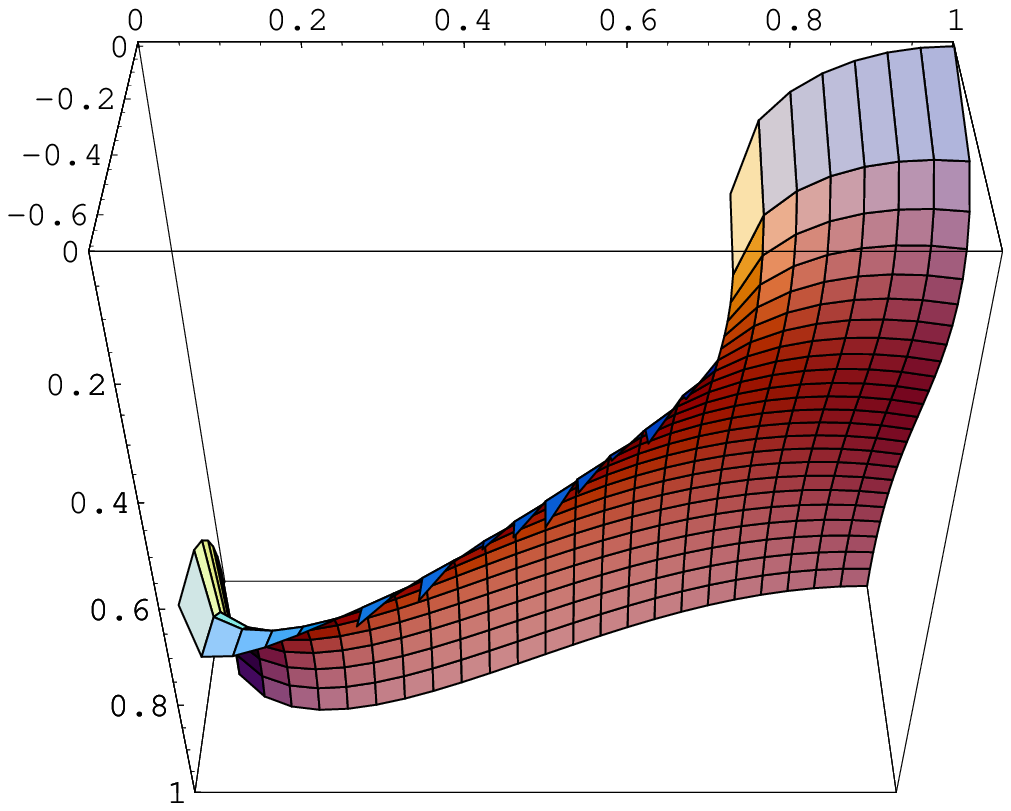}}
   \subfigure[Front-up view ]{
   \includegraphics[keepaspectratio, scale=0.35]{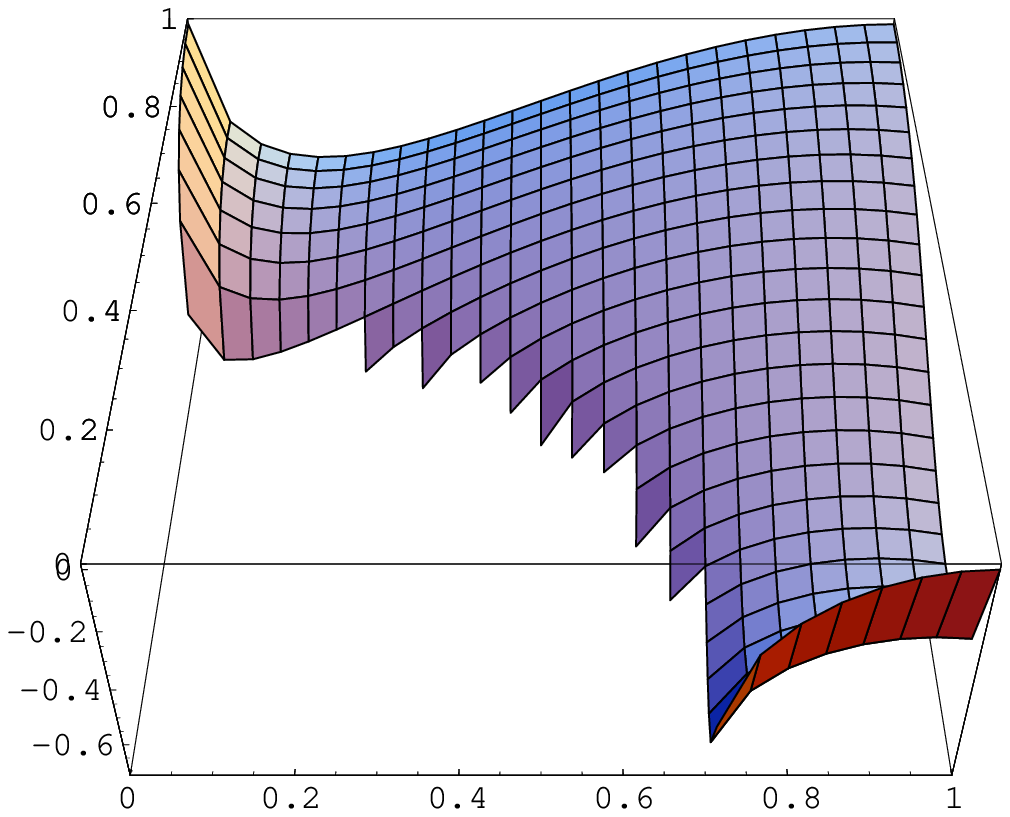}}
   \caption{Different views of the devil-fish surface}
   \end{figure}

Therefore all the cases are done and the theorem is proved.$\hfill\square$\\

\noindent The following immediate corollary is true:

\begin{Corollary} Let $a,b$ and $c$ be the sides of a triangle
$\triangle ABC$ and let also
   $m_a,m_b$ è $m_c$ be the medians to these sides
   respectively. Then
   \begin{itemize}
   \item[a)] $(2p-3a)m_a+(2p-3b)m_b+(2p-3c)m_c\ge
   0,\,\,\mbox{where}\,\,p:=\dfrac{a+b+c}{2}.$
   \item[b)]$\dfrac{m_a}{m_c}\le\dfrac{\sqrt{ab}+\sqrt{ac}+\sqrt{bc}}{a+b+c}\le
   1,\,\,a\ge b\ge c$.
   \end{itemize}
   \end{Corollary}
\bibliographystyle{plain}

\begin{thebibliography}{99}

\bibitem{Boichev}
G.~Boichev.
\newblock Geometric inequalities containing medians and some other elements of
  the triangle.
\newblock {\em Education in mathematics}, 1:41--49, 1981.
\newblock (in Bulgarian).

\bibitem{Bottema}
O.~Bottema, R.~\^Z. Djordevi\'c, R.~R. Jani\'c, D.~S. Mitrinovi\'c, and P.~M.
  Vasi\'c.
\newblock {\em Geometric inequalities}.
\newblock Groningen, 1969.

\bibitem{Nelson}
R.~B. Nelson.
\newblock Proof {Without} {Words}: The {Arithmetic-Logarithmic-Geometric}
  {Mean} {Inequality}.
\newblock {\em Math. Mag.}, 8:305, 1995.

\bibitem{Stoilov}
T.~Stoilov and K.~Chilingirova.
\newblock {\em Inequality problems}.
\newblock Narodna prosveta, Sofia, 1989.
\newblock (in Bulgarian).

\end{thebibliography}


\end{document}